\documentclass[11pt, a4paper]{article}

\usepackage[a4paper, top=2.5cm, bottom=2.5cm, left=2.5cm, right=2.5cm]{geometry}

\usepackage{amsmath}
\usepackage{booktabs}
\usepackage{enumitem}
\usepackage[sort&compress]{natbib}
\usepackage{graphicx}
\usepackage{bm}
\usepackage{caption}
\usepackage{hyperref}

\newtheorem{theorem}{Theorem}
\newtheorem{lemma}{Lemma}

\title{The Coupling Strength Is a Scale Parameter in Threshold Power-Law Reservoirs and Does Not Influence Training Accuracy}
\author{Wilten Nicola\textsuperscript{1,2,3}
\\
$^\text{1}$ Hotchkiss Brain Institute 
\\
$^\text{2}$  Department of Physics and Astronomy,\\ University of Calgary, Alberta, Canada
\\
$^\text{3}$Department of Cell Biology and Anatomy,\\ University of Calgary, Alberta, Canada}

\begin{document}
\maketitle

\begin{abstract}
In reservoir computing, the coupling strength of the initial untrained recurrent neural network (the reservoir) is an important hyperparameter that can be varied for accurate training.  A common heuristic is to set this parameter near the ``edge of chaos", where the untrained reservoir is near the transition to chaotic dynamics, and the chaos can be ``tamed".  Here, we investigate how the overall connectivity strength should be varied in threshold power-law recurrent neural networks, where the firing rate is 0 below some threshold of the current and is a power function of the current above this threshold.  These networks have been previously shown to exhibit chaotic solutions for very small coupling strengths, which may imply that the chaos cannot be tamed at all.  We show that for reservoirs constructed with threshold power-law transfer functions, if the reservoir can be trained for one single positive value of the initial reservoir coupling strength, then there exist networks with identical accuracy for all positive coupling strengths, implying that the chaotic dynamics can always be tamed or never be tamed.  This is a direct consequence of the coupling strength of threshold power-law RNNs acting as a scale parameter that does not qualitatively influence the dynamics of the system, but only scales all system solutions in magnitude.  This is independent of the power of the transfer function, with the exception of Rectified Linear Unit (ReLU) networks.  This is in contrast with conventional RNNs/reservoirs employing sigmoidal firing rates, where the strength of the recurrent coupling in the initial reservoir determines the performance on different tasks during training and also influences the network dynamics explicitly.  

\end{abstract}

\section*{Introduction} 

Recurrent Neural Networks (RNNs) serve as powerful tools in the pursuit of understanding brain operations.  When these networks are untrained and coupled randomly, they provide valuable insights into how activity propagates through neural circuits \cite{sompolinsky,ic1}, how the presence of distinct subpopulations impacts overall network dynamics \cite{yonatan}, how different time scales emerge in populations of neurons \cite{merav}, how connectivity leads to synchronizability \cite{jordan}, and the types of transitions or bifurcations neural circuits may undergo as functions of a few key parameters \cite{omri,kadmon,sompolinsky,chaosreview}.   In addition to the analysis of emergent dynamics in untrained RNNs, a parallel approach of training RNNs to mimic the behaviours of real neural circuits has also recently emerged \cite{force1,force3,palign,bptt}.  These RNNs are often trained with reservoir computing \cite{esm2,res1}, where the initially untrained RNN is set near the ``edge of chaos", a parameter regime where the network dynamics are near a critical point at the onset of self-perpetuating firing rate fluctuations (chaos).

Recently, biologically plausible threshold power-law RNNs have been developed that consider non-sigmoidal firing rates.  These networks have been studied with varying powers including square-root firing rates which mimic the firing rate typically expected from Hodgkin-Class I neurons \cite{bard}, supralinear firing rates which model responses of cells in the visual cortex \cite{ssl1,ssl2,ssl3,ssl4}, and Rectified Linear Unit (ReLU) firing rates which emerge for neurons with spike-frequency adaptation at their steady-state \cite{bardsfa}.  In parallel to the numerical simulations of these networks, variants of these networks have been analyzed from multiple perspectives including how they transition to hyperchaotic dynamics \cite{kadmon,omri,lyp1}, how slow time scales emerge in these networks \cite{bd1}, how information gets encoded by these RNNs \cite{sven}, along with analytical derivations of the moments of large-scale populations in \cite{pl2,mh}.  In the former cases, the majority of the networks were ReLU based, while supralinear firing rates were most notably analyzed in stabilized supralinear networks \cite{ssl1,ssl2,ssl3,ssl4}. 

The work in \cite{kadmon,omri} demonstrates that for large threshold power-law RNNs, the chaotic dynamics of these networks emerge immediately for any non-zero coupling strength.  Given the edge-of-chaos heuristic for conventional reservoirs, and the immediate onset of chaotic dynamics in threshold power-law RNNs, the question of whether or not threshold power-law RNNs are trainable with reservoir computing naturally emerges.  Further, if threshold power-law RNNs are trainable, how should they be initialized in terms of the coupling strength hyperparameter? 

Here, we show both analytically and numerically that threshold power-law RNNs operate differently from their classical sigmoidal cousins, but can be trained.   First, we found that the coupling strength analytically has no impact on the potential training accuracy for any supervisor, for all non-RELU threshold power-law RNNs.  This was a consequence of the fact that the coupling strength in these networks is actually a scale parameter and the network transfer functions are scale-invariant.  The coupling strength does not qualitatively alter the dynamics of the reservoir at all, beyond scaling the amplitudes of the different solutions, along with their basins of attraction.  We showed analytically and numerically that if a network could be trained for any non-ReLU threshold power-law firing rate for a single value of $g$ with some level of test accuracy, then there exists a trained reservoir for any other value of $g$ with identical accuracy.  Further, we demonstrate numerically that threshold power-law RNNs can indeed learn target dynamical systems like oscillators and low-dimensional chaotic dynamics for powers less than 1.

\section*{Results} 
\begin{figure}[htp!]
\centering
\includegraphics[scale=0.6]{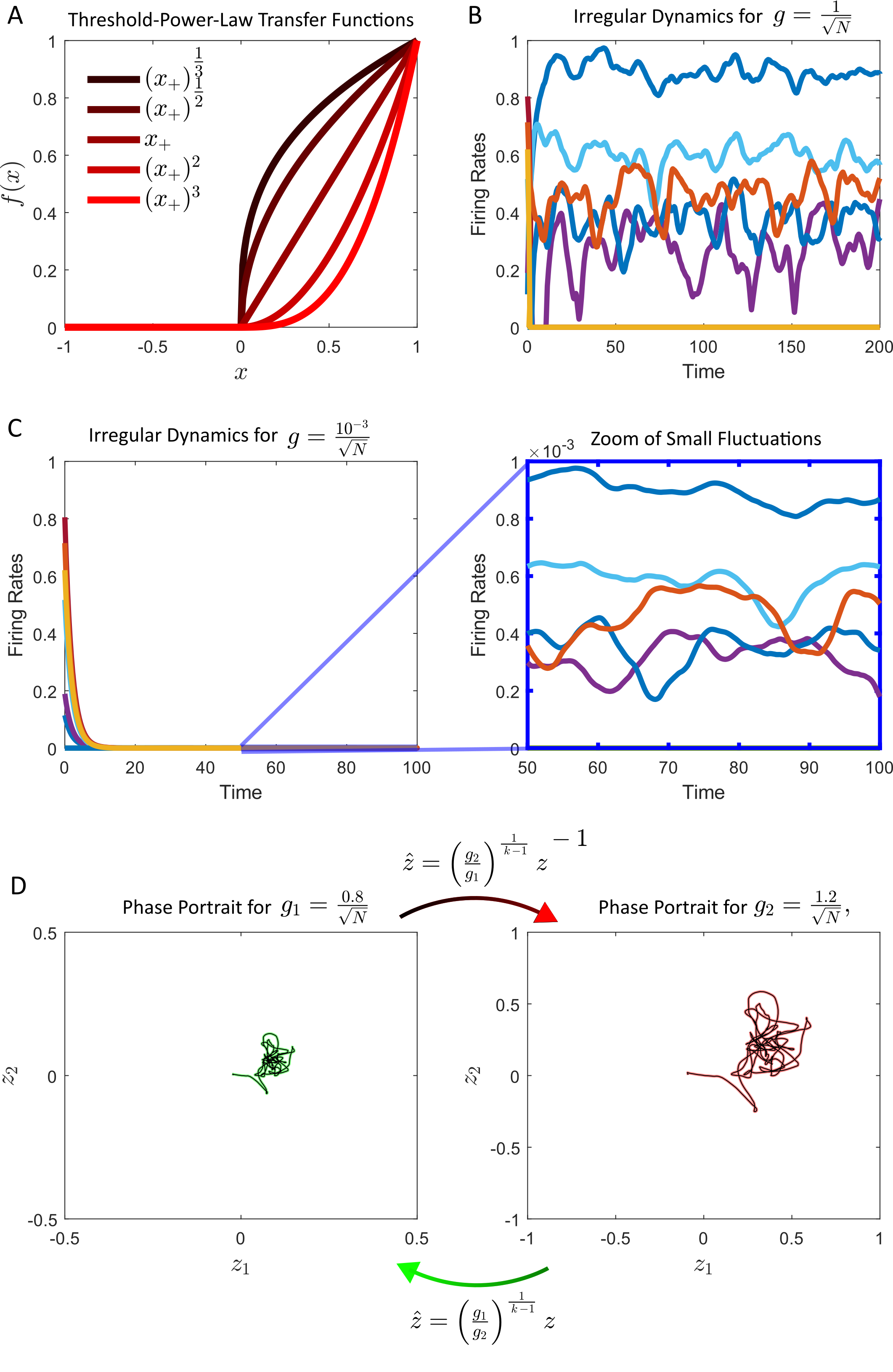}
\caption{Threshold Power-Law Recurrent Neural Networks.  \textbf{(A)} Various transfer functions for threshold power-law neural networks, where $x_+ = \max{x,0}$.  \textbf{(B)} A simulation of a threshold power-law RNN with $N=1000$ neurons, and $k=\frac{1}{2}$, where $\omega_{ij}$ is drawn from a standard normal distribution with a coupling strength of $g = \sqrt{N}^{-1}$. The system displays irregular dynamics indicative of a chaotic solution as predicted from dynamic mean field theoreis \cite{kadmon,omri} for large $N$. \textbf{(C)} Decreasing $g$ to $10^{-3}\sqrt{N}^{-1}$ leads to an initial decay to a neighbourhood near the origin (left), but a zoom (right) reveals small firing rate fluctuations.  All other parameters are as in (B).  \textbf{(D)} Lemma 1 shows that solutions for one value of $g$ are rescaled solutions for any other value of $g$.  This is shown numerically with $g_1 = 0.8\sqrt{N}^{-1}$ and $g_2 = 1.2\sqrt{N}^{-1}$ with the transform defined in the red and green arrows to rescale the respective solutions for the two values of $g$.  The initial conditions are also rescaled.  All weights were drawn from a standard normal distribution.}\label{figure1}
\end{figure}

Consider the threshold power-law RNN given by the equations 
\begin{eqnarray}
\frac{dz_i}{dt} &=& -z_i + g\sum_{j=1}^N \omega_{ij} f(z_j), \quad i=1,2,\ldots N \label{ss1} \\ 
f(z) &=& \begin{cases} z^k & z>0 \\ 0 & z \leq 0 \end{cases}, \quad k>0 \label{ss2}
\end{eqnarray}
for any $k>0$ (Figure \ref{figure1}A).   The function $f(z_j)$ is typically interpreted as the firing rate for neuron $j$ with some common examples being $k=\frac{1}{2}$ corresponding to the firing rate of the quadratic integrate-and-fire neuron \cite{bard}, while $k=1$ corresponds to the ReLU (Rectified Linear Unit) \cite{relu1} and $k>1$ corresponds to stabilized supralinear networks \cite{ssl1,ssl2,ssl3,ssl4}.  The variable $z_j$ is sometimes interpreted as the current for neuron $j$.  When $f(z)$ is a sigmoidal function and $g\propto \sqrt{N}^{-1}$, the network in (\ref{ss1}) is of the class originally considered by \cite{sompolinsky}.  For the classical rate networks with sigmoidal transfer functions considered in \cite{sompolinsky}, namely, $f(z) = \tanh(z)$, chaotic solutions emerge as $N\rightarrow \infty$ for $g>\sqrt{N}^{-1}$, while the system's equilibrium point at $z=0$ is asymptotically stable for $g<\sqrt{N}^{-1}$ and a sharp transition occurs for $g=\frac{1}{\sqrt{N}}$.  

More recently, threshold power-law recurrent networks similar to equations (\ref{ss1})-(\ref{ss2}) were considered in \cite{kadmon}. In \cite{kadmon}, the authors conclude that for powers $k$ where $k\leq \frac{1}{2}$, there exists no stable-fixed point in a system similar to (\ref{ss1})-(\ref{ss2}) as $N\rightarrow \infty$ and determine numerically that the system can be chaotic even for small values of $g$, which we reproduce here in Figure \ref{figure1}B-C.  Both \cite{kadmon} and \cite{omri} also show that as $N\rightarrow \infty$, chaos emerges even for small coupling strengths ($g\ll \sqrt{N}^{-1}$), and thus chaos is likely the only solution for large $N$, and all $g\neq 0$.   

We considered the non-sigmoidal firing rates in (\ref{ss1})-(\ref{ss2}), which were more generally analyzed in \cite{kadmon,omri}, for their potential use in reservoir computing.   In many instantiations of reservoir computing, the weights of (\ref{ss1}) are modified with a rank-$m$ perturbation to stabilize the system onto the dynamics of some $m$-dimensional supervisor, $\hat{\bm x}(t)$
\begin{eqnarray}
 \frac{dz_i}{dt} &=& -z_i + g\sum_{j=1}^N \omega_{ij} f(z_j)+ 
\sum_{l=1}^m\bm\eta_{lj} \bm \phi_{lj} f(z_j)\\
\hat{\bm x}(t) &=&  \bm \phi^T f(\bm z) 
\end{eqnarray}
where $\hat{\bm x}(t)$ is a $m$ dimensional approximation to the target dynamics of $\bm x(t)$ and the encoders ($\bm \eta_i$) and decoders ($\bm \phi_i$) are $m\times 1$ vectors for each neuron (with index $i$).  The initial weights, $\omega_{ij}$ are randomly generated from the Gaussian with mean 0, and standard deviation 1.  The $\frac{1}{\sqrt{N}}$ that is conventionally present can be absorbed into the scalar $g$.  Note that negative values of $g$ typically need not be considered, as the sign of $g$ can be absorbed into the weights $\omega_{ij}$.  In the case of weights drawn from random distributions with mean 0, as in \cite{omri,kadmon}, the dynamics of the network are unaffected.   

In matrix form, the reservoir is the initial nonlinear dynamical system given by the equation:
\begin{eqnarray}
\frac{d\bm z}{dt} &=&-\bm z + g\bm\omega f(\bm z) \quad \text{(Initial Reservoir)} \label{r1}
\end{eqnarray}
where after successful training, the system's weights are perturbed with the low-rank matrix $\bm \eta \bm \phi^T$: 
\begin{eqnarray}
\frac{d\bm z}{dt} &=& -\bm z + \left(g\bm \omega +\bm \eta\bm\phi^T \right)f(\bm z), \quad \hat{\bm  x} = \bm \phi^T f(\bm z) \quad \text{(Trained Reservoir)} \label{r2}
\end{eqnarray}
and $\hat{\bm x}(t) \approx \bm x(t)$.  During learning, $\bm \eta$ is typically held fixed and is randomly generated at the start of a simulation, while $\bm \phi$ is learned in some way to minimize the following loss function (the mean-squared error):
\begin{eqnarray}
 L(\bm \phi,\bm \eta,g,\bm \omega) = \sum_{l=1}^m\int_0^T \left(\hat{x}_l(t') - x_l(t')\right)^2\,dt' \label{test_loss}
\end{eqnarray}
although other losses may also be considered.  Note that we've explicitly identified all parameters $(\bm \phi,\bm \eta,g,\bm \omega)$ which determine the loss. We will interpret the loss in equation (\ref{test_loss}) as the test loss where the interval $[0,T]$ occurs after training has concluded (and $\bm \phi$ is held constant).  There are many techniques to learn the decoder $\bm \phi$  either iteratively in an online setting (e.g. FORCE training \cite{force1,force2,force3} or Backprop through time \cite{bptt}) or offline (e.g. Echo State Training \cite{esm1,esm2}), or more recently, through predictive alignment \cite{palign}.  We do not specify how the decoder is learned in the analysis, but we use recursive least squares (RLS) here (Materials and Methods, \cite{force1,force2}). 

When $f(z)=\tanh(z)$, typically $g$ is selected near $g=\frac{1}{\sqrt{N}}$ as this is the ``edge of chaos" for large networks, although the range of $g$ required varies with the learning algorithm employed and the supervisor considered \cite{force1}. When the $\frac{1}{\sqrt{N}}$ term is absorbed into the weights, the edge of chaos occurs for $g=1$ \cite{sompolinsky}.   
In this work, we seek to determine how $g$ should be scaled for threshold power-law reservoirs, given previous reports of chaos for even small values of $g$ for large networks \cite{kadmon}.  Surprisingly, we find that $g$ does not explicitly determine the loss, for any-sized network, as stated in the following theorem:  
\begin{theorem}

Suppose that the initial reservoir (\ref{r1}) is trained on some $m$-dimensional supervisor $\bm x(t)$ with some encoder/decoder pair $\bm \eta$, $\bm \phi$ for $g = g^*$, achieving a test-loss of $ L^* = L(\bm \phi,\bm \eta,g^*,\bm \omega) $.  Further, assume that the threshold power-law considered is $k\neq 1$ (non-ReLU).  Then for any $\bar{g}>0$, there exists an encoder and decoder pair, $\hat{\bm \eta}$, $\hat{\bm \phi}$ with $L(\hat{\bm \phi},\hat{\bm \eta},\bar{g},\bm \omega)= L^*$ 
\end{theorem}

Before we prove Theorem 1, there are a few key points to note. First, if the initial reservoir can be trained with some $g=g^*$, then it can be trained for all $g$ with equivalent accuracy $L = L^*$, so long as the correct encoder and decoder pair ($\hat{\bm \eta}$, $\hat{\bm \phi}$) can be found by the learning algorithm employed.  Second, through the very straightforward derivation of Theorem 1, we will demonstrate that $\hat{\bm \eta}$, $\hat{\bm \phi}$ can be expressed in terms of $g,g^*,{\bm \eta}$, ${\bm \phi}$, and $k$.  The third key point to note is that Theorem 1 implies that if a learning algorithm can find an optimum (local or global) low-rank perturbation matrix for a single value of the coupling strength $g$ in (\ref{r2}), then this solution is also a local/global optimum for any $g$.  Finally, the theorem does not specify how the initial reservoir weights $\bm \omega$ are generated, or even if the training was successful ($L_0$ is small).  Thus, if the chaotic dynamics can be tamed for a single value of $g=g^*$, they can be tamed for all values of $g$.   The reason this theorem can be so general is that it is a consequence of the following lemma:

\begin{lemma}
Consider the threshold power-law recurrent neural network given by 
\begin{eqnarray}
\frac{d\bm z}{dt} = -\bm z + g\bm \omega f(\bm z) \label{sys1}
\end{eqnarray}
where $k\neq 1$.  Then (\ref{sys1}) can be rescaled to the following system:
\begin{eqnarray}
\frac{d\bm y}{dt} = -\bm y + \bm \omega f(\bm y) \label{sys1}
\end{eqnarray}
with $\bm y = g^{\frac{1}{k-1}}\bm z$ where $k$ is the power of the threshold power-law transfer function. 
\end{lemma}
The derivation of Lemma 1 is extremely straightforward: 
\begin{eqnarray*}
\frac{d\bm y}{dt} &=& g^{(k-1)^{-1}} \frac{d\bm z}{dt} \\
&=&  g^{(k-1)^{-1}} (\bm y g^{-(k-1)^{-1}} + g\bm \omega f(\bm y g{-(k-1)^{-1}})   \\
&=& -\bm y + g^{(k-1)^{-1} + 1 -k(k-1)^{-1}}\bm \omega f(\bm y) \\
&=& -\bm y +  \bm\omega  f(\bm y)
\end{eqnarray*}
In plain terms, Lemma 1 states that the coupling strength, $g$, acts as a scale parameter for threshold power RNNs when $k\neq 1$. Thus, $g$ does not qualitatively influence the network dynamics for all $g\neq 0$, and only scales the system's solutions.  All solutions that exist for a particular $\bm \omega$ at some particular value of $g$ exist for all values of $g$ for that same $\omega$, up to a rescaling factor.  This has multiple implications, including the straightforward derivation of Theorem 1.  First, if the system has chaotic solutions for a particular value of $g^*$, then it has (the same) chaotic solutions for all $g\neq g^*$ and these are rescaled versions of the chaotic solutions that exist for $g=g^*$.  This was numerically confirmed in simulations (Figure \ref{figure1}D).   There is no transition to chaos by using the connectivity strength of the weights $g$ as a bifurcation parameter in non-ReLU threshold power-law RNNs.  For any finite $N$ and $\bm \omega$, if the system displays chaotic dynamics for some $g$, it does so for all $g>0$.  Indeed, in \cite{kadmon}, the authors remark that the chaos found in the threshold-power-law networks for $k=\frac{1}{2}$ can persist for very small $g$ (see Figure 2a in \cite{kadmon}) although we remark that the network in \cite{kadmon} differs from the RNNs considered here as the authors consider external inputs rather than completely autonomous dynamics.  Despite the (very) simple derivation of Lemma 1, it does not appear to be explicitly mentioned in the literature surrounding threshold power-law RNNs to the best of our knowledge. 

\begin{figure}[htp!]
\centering
\includegraphics[scale=0.8]{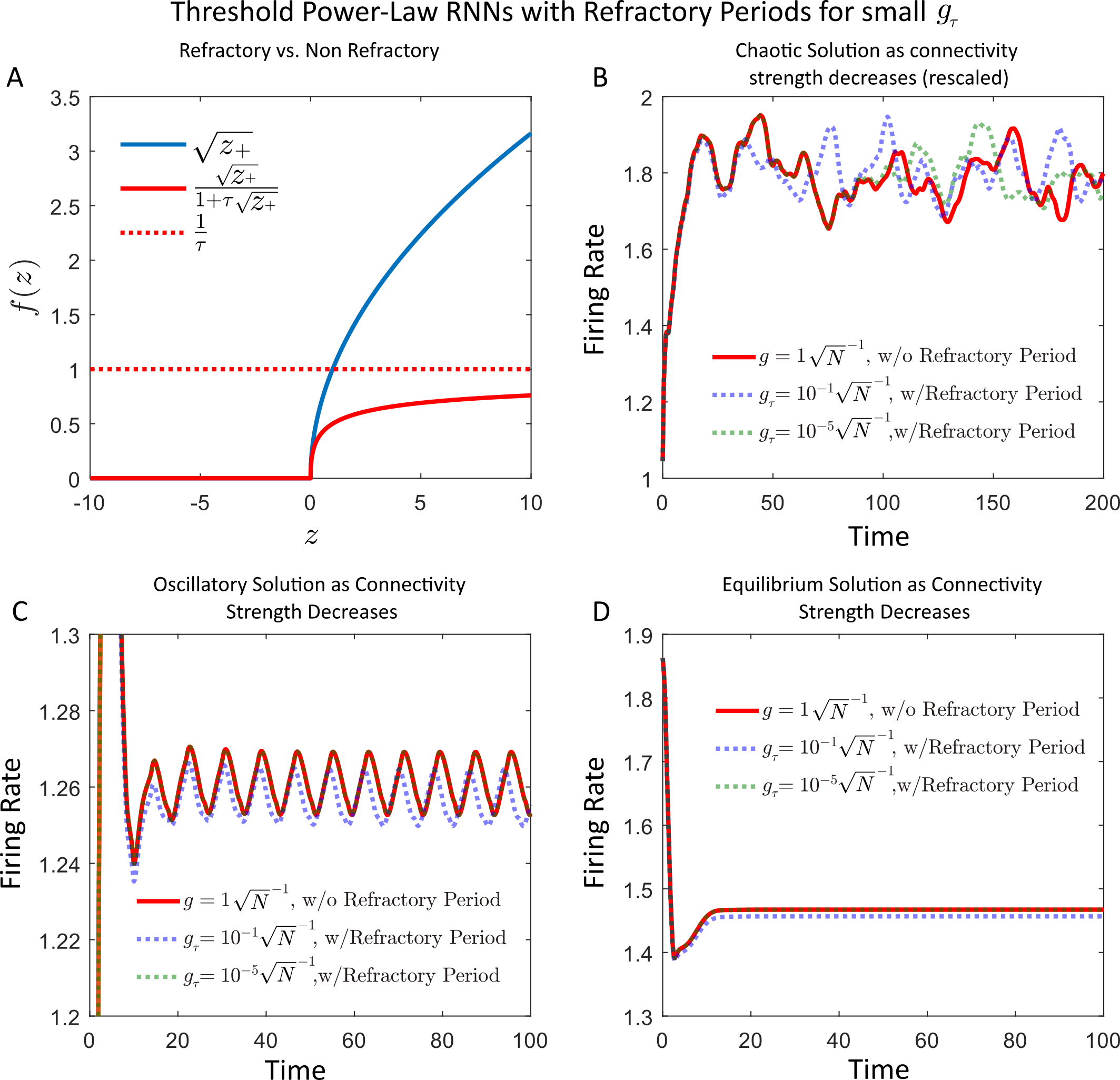}
\caption{Threshold Power-Law RNNs with Refractory Periods.  \textbf{(A)} The firing rates for a threshold power-law RNN with $k=\frac{1}{2}$ (blue) and with a refractory period (red).  The firing with a refractory period $\tau$ converges to $\tau^{-1}$ for high input currents ($z$).  \textbf{(B)} The firing rate for a simulated network without a refractory period (red, $N=2000$, $g = \sqrt{N}^{-1}$ and with a refractory period for $g_\tau =10^{-1}\sqrt{N}^{-1}$ (blue-dashed) and $g_\tau=10^{-5} \sqrt{N}^{-1}$  (green-dashed).  As $g_\tau\rightarrow 0$, the solutions of the network with a refractory period converge to solutions of the network without one. Note that the green and blue solutions were re-scaled with $\hat{z}(t) =  z(t)(\frac{g_{\tau}}{g})^{1/(k-1)}$ where $g_\tau$ denotes the network coupling strength for the network with a refractory period $\tau$. \textbf{(C)} Identical as in (B) only with oscillatory solutions as $g_\tau\rightarrow 0$ by with $N=50$ neurons.  \textbf{(D)} Identical as in (C), only with an equilibrium point solution as $g_\tau \rightarrow 0$.  All weights were drawn from a standard normal distribution. }\label{figure2}
\end{figure}

Further, we remark that Lemma 1 also has implications for reservoirs where the firing rates saturate due to refractory periods \cite{bard}.  The firing rate for these networks is given by 
\begin{eqnarray}
f_\tau(z) = \frac{f(z)}{\tau f(z) + 1 } \label{ref4}
\end{eqnarray}
which asymptotically converges to $f_\tau(z)\rightarrow \frac{1}{\tau}$ as $z\rightarrow \infty$ and $f_\tau (z) \rightarrow f(z)$ as $z\rightarrow 0$.  Consider the threshold power-law RNN with a refractory period given by 
\begin{eqnarray}
\frac{d \bm z}{dt} = -\bm z + g_\tau \bm \omega f_\tau (\bm z)  \label{ref1}
\end{eqnarray}
coupled with coupling parameter $g_\tau$, where $f(z)$ is a threshold power-law transfer function as in Equation (\ref{ss2}) while $f_\tau(z)$ is given by Equation (\ref{ref4}).  In the limit that $g_\tau \rightarrow 0$, the solution(s) to equation (\ref{ref1}) can be shown to converge to rescaled solutions of
\begin{eqnarray}
\frac{d\bm y}{dt} = -\bm y + \bm \omega f(\bm y).
\end{eqnarray}
This was confirmed numerically in Figure \ref{figure2} for chaotic (Figure \ref{figure2}B), oscillatory (Figure \ref{figure2}C), and equilibria solutions (Figure \ref{figure2}D).

The proof of Theorem 1 follows immediately from Lemma 1 by considering the trained reservoir:
\begin{eqnarray}
\frac{d\bm z}{dt} = -\bm z + (g^*\bm \omega + \bm \eta \bm \phi^T) f(\bm z), \quad \hat{\bm x}  = \bm \phi^T \bm f(\bm z)
\end{eqnarray}
which is equivalent to the rescaled system 
\begin{eqnarray}
\frac{d\bm y}{dt} = -\bm y + (g\bm \omega + \hat{\bm \eta} \hat{\bm \phi}^T) f(\bm y), \quad  \hat{\bm x}  = \bm \hat{\bm\phi}^T  f(\bm y)
\end{eqnarray}
where 
\begin{eqnarray}
\hat{\bm \phi} = \bm \phi \left(\frac{g}{g^*}\right)^{\frac{k}{k-1}},\quad \hat{\bm \eta} = \bm \eta \left(\frac{g^*}{g}\right)^{\frac{1}{k-1}} \label{tform}
\end{eqnarray}

Thus, if a threshold power-law reservoir can be trained for a single coupling strength $g$ (for a particular supervisor), it can be trained for all coupling strengths. 

\begin{figure}[htp!]
\centering
\includegraphics[scale=0.39]{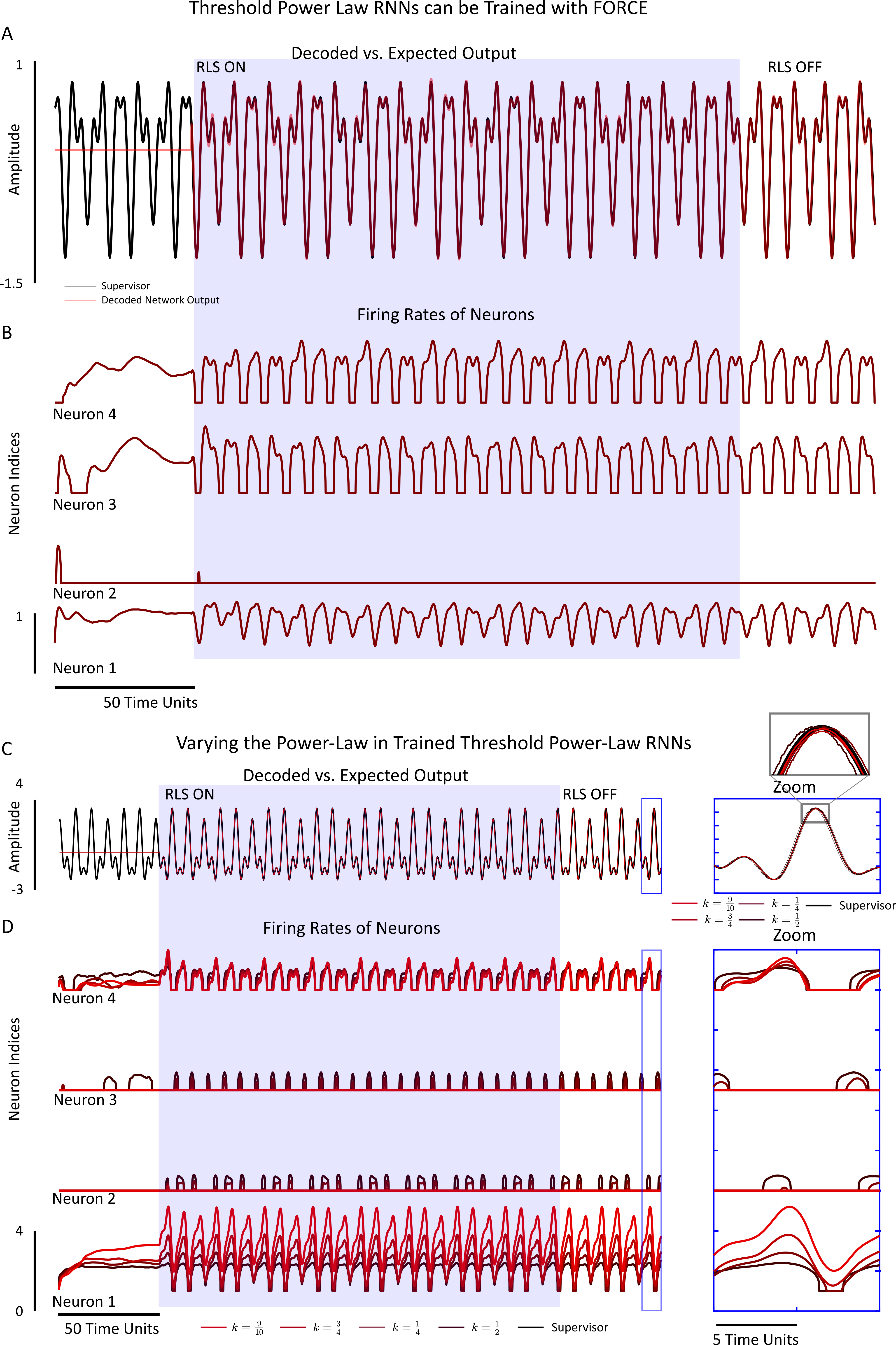}
\caption{Threshold Power-Law RNNs can be Trained with FORCE.  \textbf{(A)} The supervisor (black), a random sum of oscillators (Materials and Methods) vs the network approximation (red).  RLS was turned on after 50 time units and turned off for the last 50 time units in a 300 time unit simulation.  The network consisted of $N=2000$ neurons with $k=\frac{1}{2}$ and $g = 1$.  \textbf{(B)} The firing rates for 4 neurons, time aligned with (A).  \textbf{(C)} The supervisor (as in (A), but with a different random sum, black) vs networks with increasingl larger powers $k$ (reds).  RLS was turned on after 50 time units and off for the last 50 time units in a 300 time unit simulation.  A zoom of the last 10 time units is shown on the right.  \textbf{(D)} The firing rates for four neurons, time-aligned with (C).  Note that all 4 networks used the same initial state and the same initial reservoir weight matrix with $g=1.5\sqrt{N}^{-1}$ in all cases.  A zoom of the last 10 time units is shown on the right.  All weights were drawn from a standard normal distribution. 
}\label{figure3}
\end{figure}

\begin{figure}[htp!]
\centering
\includegraphics[scale=0.8]{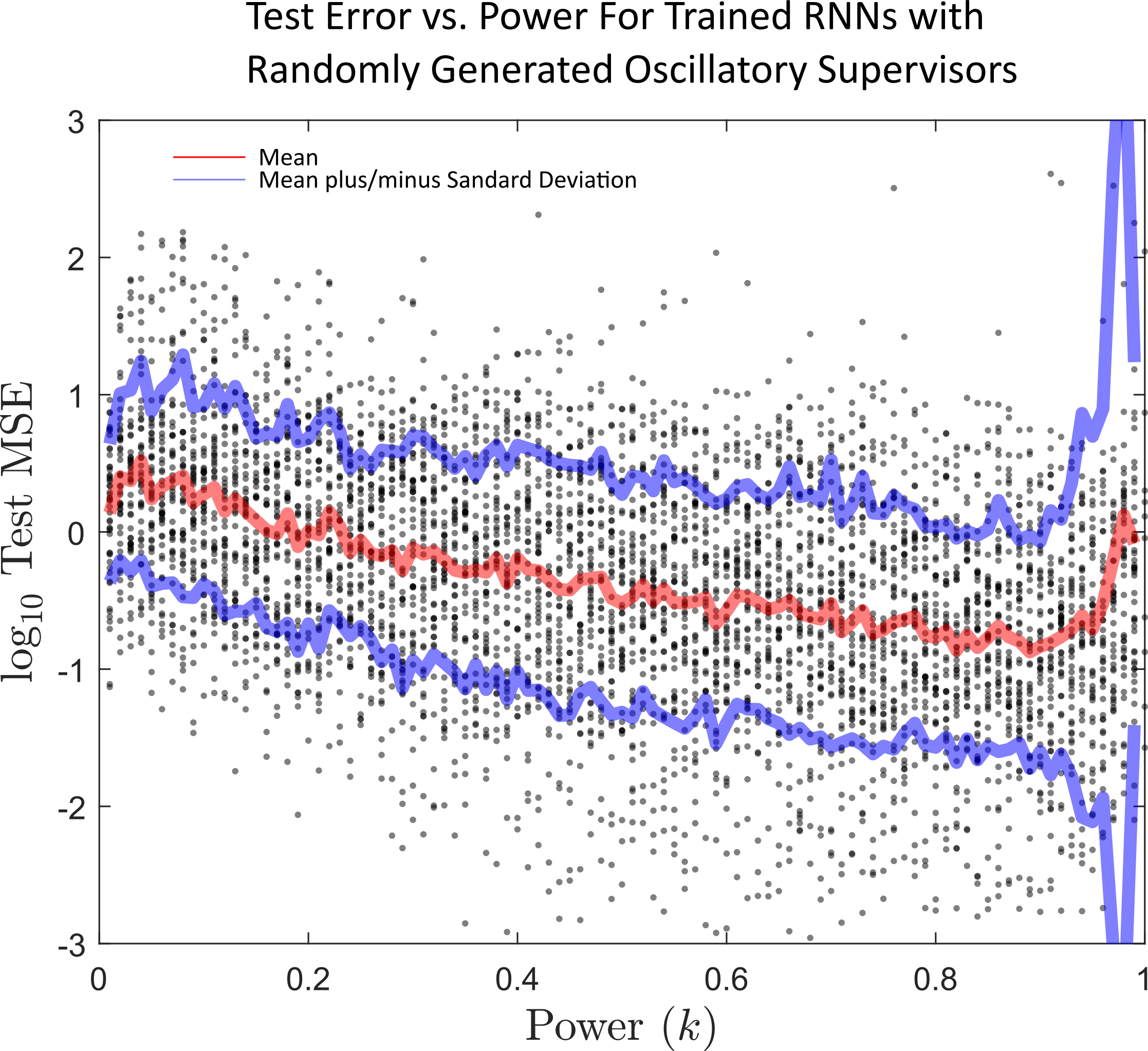}
\caption{Test Error vs. Power in Trained Threshold Power-Law RNNs.  The power $k$ varied from $1/100$ to $100/100$ in increments of $1/100$.  For each value of $k$, 50 random networks consisting of $N=2000$ neurons were generated and FORCE trained with randomly generated oscillatory supervisors as in Figure \ref{figure3}.  The log of the test MSE is shown as individual points (black dots) with the mean (red) and mean $\pm$ standard deviation in blue.  The test error monotonically decreases until $k\approx 1$ where the networks start becoming unstable.  Each network was simulated for 500 time units with RLS turned on after the first 50 time units and off for the last 50 time units (testing). All networks used a constant $g =1.5 \sqrt{N}^{-1}$.   All weights were drawn from a standard normal distribution.  }\label{figure4}
\end{figure}

Next, we sought to determine if these threshold power-law RNNs could be trained at all, that is if the chaos could be tamed onto a target dynamical system.   We found that threshold power-law RNNs (for $k=1/2$) could indeed learn different dynamical systems with FORCE training (Figure \ref{figure3}) \cite{force1}.  We first considered training the networks to mimic multi-frequency oscillations produced from randomly generated oscillators of the form 
\begin{eqnarray} 
x(t) = \sum_{j=1}^{n_o} a_j \cos(2\pi t f_j)
\end{eqnarray} 
where $a_j$ is a randomly generated amplitude drawn from a standard normal, and $f_j$ is a fixed frequency (Materials and Methods). We found that networks could be trained for $k=\frac{1}{2}$ (Figure \ref{figure3}A-B), and with other powers for $k<1$ (Figure \ref{figure3}C-D).   We found that using $k>1$ led either to unstable system dynamics (blow-up, not shown) or convergence to non-firing dynamics (not shown).   Further, we found that the test error in training randomly generated oscillatory supervisors decreased steadily with power $k$ (Figure \ref{figure4}) until $k\rightarrow 1$, where the reservoirs could lose stability or testing accuracy.  Thus, threshold power-law RNNs can be trained with FORCE reliably for sufficiently sublinear powers for the randomly generated oscillatory supervisors considered here.  

\begin{figure}[htp!]
\centering
\includegraphics[scale=0.5]{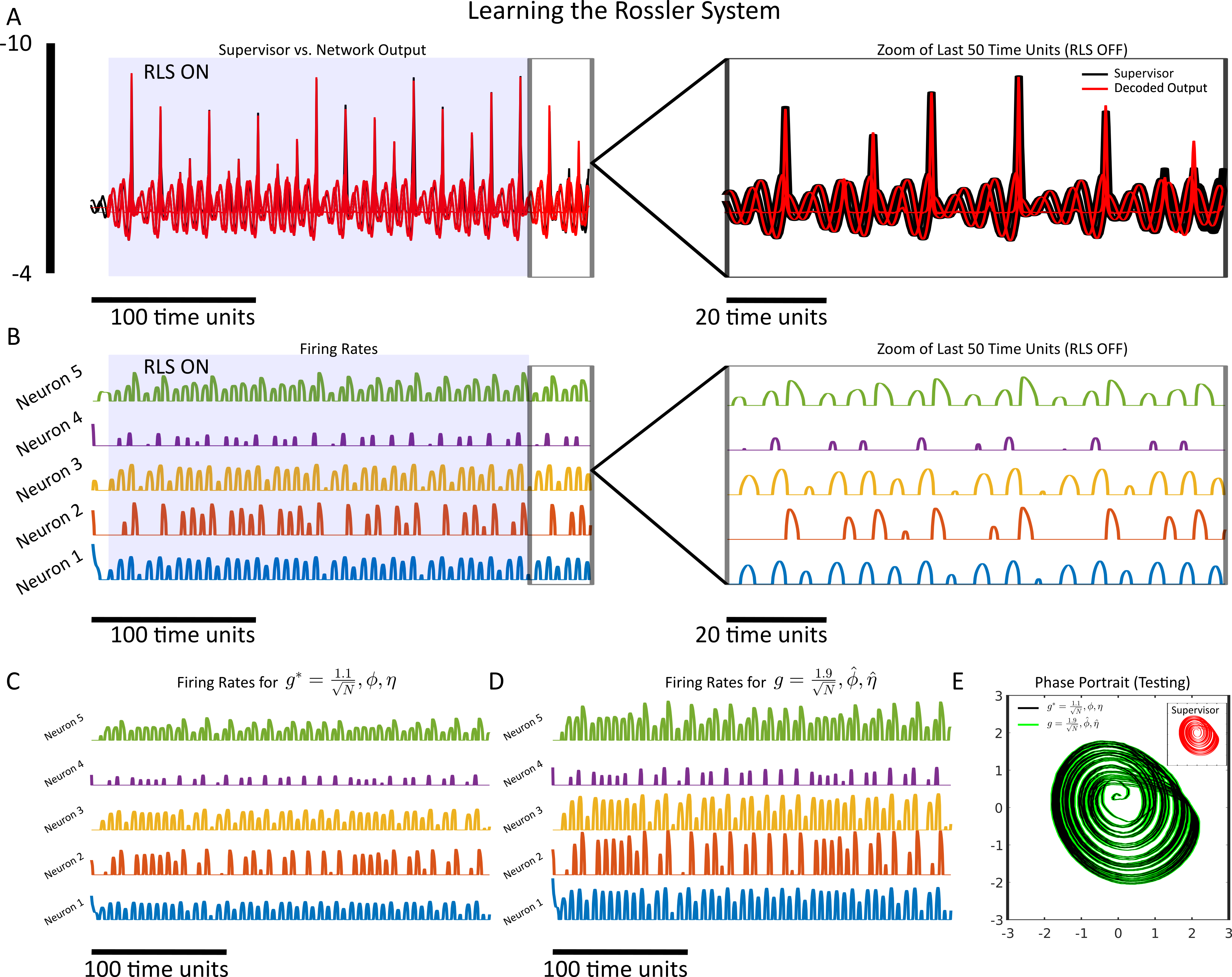}
\caption{Trained Threshold-Power-Law Recurrent Neural Networks can be Rescaled to Any Reservoir Strength.  \textbf{(A)} A network of $2000$ neurons with $g = \frac{1.1}{\sqrt{N}}$ was trained on a complex oscillator task with FORCE training.  RLS was turned on on $50$ time units and turned off at $200$ time units. The decoded output (red) vs the target supervisor (black) is shown on top, while the firing rates for 5 neurons are shown on the bottom.  \textbf{(B)} A network of $2000$ neurons with $g = \frac{1.1}{\sqrt{N}}$ was trained to learn the Rossler dynamical system with FORCE training (network, red, system, Rossler, black).  A zoom of the last 50 time units with RLS off is shown on the right, while the full simulation is shown on on the left.  \textbf{(C)} The firing rates for 5 randomly selected neuron.  The zoom on the right is aligned with the zoom in (B). \textbf{(D)} The firing rates for the reservoir with $g^* = \frac{1.1}{\sqrt{N}}$ during testing.  \textbf{(E)} The firing rates for the reservoir with $g = \frac{1.0}{\sqrt{N}}$ during testing using the encoder and decoder pair $\hat{\bm\eta}$ and $\hat{\bm \phi}$ where $\phi$ as described in the text.  \textbf{(F)} The phase portrait of the $x_1$ vs $x_2$ variable for the network trained with $g^*$ (black) and $g$ (green).  The Rossler inset is shown in red.      }\label{figure5}
\end{figure}

Next, we tested threshold power-law reservoirs on accurate chaotic time series predictions (Figure \ref{figure5}(B-F) with the Rossler system.  We successfully trained a threshold power-law RNN with FORCE training to mimic the Rossler system (Figure \ref{figure5}A-F) with a short period of training.  Further, as confirmed by theorem 1, the trained reservoir for one value of the coupling strength $(g^* = \frac{1.1}{\sqrt{N}})$ could be rescaled to a different value of the coupling strength ($g =\frac{1.9}{\sqrt{N}}$) with equation (\ref{tform}) rescaling the encoders and decoders.  

Collectively, these results analytically and numerically confirm that the coupling strength $g$ acts as a scale parameter for threshold power-law RNNs.  This implies that when used as a reservoir, the coupling strength, $g$, does not dictate the trainability of these types of reservoirs for any supervisor.

\section*{Discussion}

Through both analytical derivations and numerical simulations, we have shown that threshold power-law RNNs, both trained and untrained, operate differently from their sigmoidal counterparts \cite{sompolinsky}.  In particular, the overall coupling strength of these networks does not influence the chaotic dynamics of the system, and as a result, does not influence the accuracy of any trained network.  The coupling strength scales solutions up or down in amplitude.  The scaling invariance in threshold power-law RNNs even applies to networks where the firing rate has a refractory period, as all weak coupling solutions $g_\tau \rightarrow 0$ are scaled solutions of the threshold power-law RNN where the firing rates have no refractory periods. 

We remark that the fact that $g$ acts as a scale-parameter was partially considered in \cite{omri} and \cite{kadmon}.  In particular, both studies show that for large network limits, chaotic solutions emerge for all $g\neq 0$.  Here, we show that this result is more general: for any finite $N$, the solutions that exist for some $g=g^*$ exist for all $g>0$, chaotic or otherwise, as the coupling strength does not influence the qualitative behaviours of threshold power-law RNNs.  Similarly, if the non-firing equilibrium point ($\bm z=0$) is unstable for any finite $N$ for the weight matrix $g^* \bm\omega$, then it is unstable for all coupling strengths with the same $\bm \omega$.  The dynamic mean-field theories considered by \cite{omri} and \cite{kadmon} guarantee the existence of high-dimensional chaotic solutions based on excitatory/inhibitory balance, similar to those of classical sigmoidal networks in \cite{sompolinsky} as $N\rightarrow \infty$.  The work here shows that the immediate onset to chaos does not impact the ability of these networks to be trained with the techniques of reservoir computing, and that the asymptotic derivations for large $N$ are caused by the coupling strength acting as a scale parameter.  We remark that threshold power-law neural networks have also been considered previously for certain values of $k$ under other conditions \cite{relu_rnn_1,relu_rnn_2,square_1,nicola2016}. 

While we have demonstrated that the parameter $g$ does not directly influence the loss of a trained reservoir or the dynamics of an untrained initial reservoir, this does not imply that the $g$ parameter is not important in the context of explicitly training an RNN.  First, we remark that Theorem 1 states that the learned low-rank perturbation for $g=g^*$ can be transformed so that the initial reservoir had a coupling strength of $g$, but this does not imply that a learning algorithm will necessarily learn the appropriate scaling for an encoder/decoder pair.  While $g$ does not impact the dynamics of the reservoir qualitatively, this scaling is still practically important for training.  If $g$ is too small, and the learning algorithm applied does not somehow explicitly scale the encoder/decoder pair to account for $g$, then a functional low-rank perturbation that stabilizes the dynamics of the reservoir to $\hat{\bm x}\approx \bm x(t)$ may not be found. 

Collectively, this work demonstrates that threshold power-law RNNS, which in some ways are more biologically plausible than sigmoidal RNNs, have very different qualitative behaviours.  The chaos that emerges for large networks can either always be tamed for any coupling strength, or can never be tamed. 

\section*{Acknowledgments} 
We thank Claudia Clopath for her insightful comments.  This work was funded by an NSERC Discovery Grant (RGPIN/04568-2020), a Canada Research Chair (CRC-2019-00416), and the Hotchkiss Brain Institute. 


\clearpage
\section*{Materials and Methods} 

\subsection*{FORCE Training}
Chaotic RNNs were trained with the FORCE algorithm \cite{force1,force2,force3}.  In FORCE training, an initially chaotic RNN given by 
\begin{eqnarray}
    \frac{d\bm z}{dt} = -\bm z + g\bm \omega f(\bm z) 
\end{eqnarray}
is stabilized onto an $m$-dimensional supervisor $\bm x(t)$ with a linear decoder, $\bm \phi$ where $\hat{\bm x(t)} = \bm \phi^T f(\bm z(t))$, and $\hat{x}(t)\approx x(t)$.  After training, the system equations are given by:  
\begin{eqnarray}
    \frac{d\bm z}{dt} = -\bm z + (g\bm \omega+\bm \eta \bm \phi) f(\bm z) 
\end{eqnarray}
The encoder matrix was an $N \times k$ matrix where each element was drawn from a uniform random distribution defined on $[-1,1]$.  The decoder was updated with Recursive Least Squares (RLS) 
\begin{eqnarray}
\bm P_{n+1} &=& \bm P_n - \frac{\bm r_n ^T \bm P_n ^T \bm P_n \bm r_n}{1+\bm r_n^T \bm P_n \bm r_n}\\
\bm \phi_{n+1} &=& \bm \phi_n - \frac{\bm P_n \bm e_n^T}{1+\bm r_n^T \bm P_n \bm r_n}
\end{eqnarray}
where $e_n = \hat{x}_n - x_n$.  The index $n$ corresponds to every 3 time steps in the numerical integration here.  All networks considered were integrated with a forward Euler scheme with a time step of $\Delta = 10^{-2}$ in MATLAB 2023a.   \subsection*{Supervisors}

The random oscillatory supervisors (Figure \ref{figure2}, \ref{figure3} were generated as 
\begin{eqnarray}
x(t) = \sum_{j=1}^{3} a_j \cos(2\pi t f_j)
\end{eqnarray}
where $f_1 = 1/6$, $f_2 = 1/8$ and $f_3 = 1/10$.  The Rossler system \cite{rossler} was given by
\begin{eqnarray}
\frac{dx_1(t)}{dt} &=& -x_2(t) - x_3(t)\\
\frac{dx_2(t)}{dt} &=& x_1(t) + ax_2(t)\\
\frac{dx_3(t)}{dt} &=& b + x_3(t)(x_1(t)-c)
\end{eqnarray}
where $a=0.2,b=0.2,c=5.7$ and $x_1(t),x_2(t)$ and $x_3(t)$ were $z$-transformed (mean 0, standard deviation 1) before being used as a supervisor.

\clearpage

\section*{Figures}

\clearpage

\clearpage

\bibliographystyle{unsrt}
\bibliography{References}

\end{document}